\documentclass[12pt,A4paper]{article}
\usepackage{amsmath,amsfonts,amsthm}

\newtheorem{thm}{Theorem}[section]

\newtheorem{cor}[thm]{Corollary}
\def\pf{\noindent{\it Proof.} }
\newcommand{\bfa}{\mathbf{a}}
\newcommand{\bfb}{\mathbf{b}}
\newcommand{\bfc}{\mathbf{c}}
\newcommand{\bfi}{\mathbf{i}}
\newcommand{\bfj}{\mathbf{j}}
\newcommand{\bfx}{\mathbf{x}}
\newcommand{\bfy}{\mathbf{y}}
\newcommand{\bfw}{\mathbf{w}}

\def\qed{\nopagebreak\hfill{\rule{4pt}{7pt}}\medbreak}

\makeatletter \@addtoreset{equation}{section} \makeatother

\begin{document}
\begin{center}
{\Large\bf Minimally Intersecting Set Partitions of Type $B$}
\end{center}

\begin{center}
William Y.C. Chen and David G.L. Wang\\[5pt]
Center for Combinatorics, LPMC-TJKLC\\[5pt]
Nankai University, Tianjin 300071, P.R. China\\[5pt]
chen@nankai.edu.cn, wgl@cfc.nankai.edu.cn
\end{center}

\begin{abstract}
Motivated by Pittel's study of minimally intersecting set
partitions, we investigate minimally intersecting set partitions of
type $B$. We find a formula for the number of minimally intersecting
$r$-tuples of $B_n$-partitions, as well as a formula for the number
of minimally intersecting $r$-tuples of $B_n$-partitions without
zero-block. As a consequence, it follows the formula of Benoumhani
for the Dowling number in analogy to Dobi\'{n}ski's formula.
\end{abstract}



\section{Introduction}

This paper is primarily concerned with the meet structure of the
lattice of type $B_n$ partitions of the set $[\pm n]=\{\pm 1,\,\pm
2,\,\ldots,\,\pm n\}$, as well as of the meet-semilattice of type
$B_n$ partitions without zero-block. The lattice structure of type
$B_n$ set partitions has been studied by Reiner \cite{Rei97}. It can
be regarded as a representation of the intersection lattice of the
type $B$ Coxeter arrangements, see Bj\"{o}rner and
Wachs~\cite{BW04}, Bj\"{o}rner and Brenti~\cite{BB05} and
Humphreys~\cite{Hum90}.

We establish a formula for the number of $B_n$-partitions $\pi'$
which minimally intersect a given $B_n$-partition $\pi$. Using the
same technique, we derive a formula for the number of
$B_n$-partitions $\pi'$ without zero-block which minimally intersect
a given $B_n$-partition $\pi$ without zero-block. The ordinary case
has been studied by Pittel~\cite{Pit00}. In particular, if we take
$\pi$ to be the minimal $B_n$-partition, our formula reduces to a
formula of Benoumhani~\cite{Ben96} for the number of
$B_n$-partitions (called the Dowling number, see
Dowling~\cite{Dow73}),   which is analogous to Dobi\'nski's formula
for the number of partitions of a finite set.

In a more general setting, we derive two formulas for the number of
minimally intersecting $r$-tuples of $B_n$-partitions and the number
of minimally intersecting $r$-tuples of $B_n$-partitions without
zero-block. Recall that Canfield~\cite{Can01} has found a relation
between the exponential generating function of the number of
minimally intersecting $r$-tuples of partitions and the powers of
the Bell numbers. We give a type $B$ analogue of this relation.

Let us give an overview of relevant notation and terminology. A
partition of a set $S$ is a collection $\{B_1,B_2,\ldots,B_k\}$ of
subsets of $S$ such that $B_1\cup B_2\cup\cdots\cup B_k=S$ and
$B_i\cap B_j=\emptyset$ for any $i\ne j$. A {\em set partition of
type $B_n$} is a partition $\pi$ of the set $[\pm n]$ into blocks
satisfying the following conditions:
\begin{itemize}
\item[(i)]%
For any block $B$ of $\pi$, its opposite $-B$ obtained by negating
all elements of $B$ is also a block of $\pi$;
\item[(ii)]%
There is at most one {\em zero-block}, which is defined to be a
block $B$ such that $B=-B$.
\end{itemize}
We call $\pm B$ a {\em block pair} of $\pi$ if  $B$ is a non-zero-block of $\pi$. For example,
\begin{align*}
\pi_1&=%
\left\{\{\pm 1,\,\pm 2,\,\pm 5,\,\pm 8,\,\pm 12\},%
\ \pm\{3,\,11\},\ \pm\{4,\, -7,\,9,\,10\},\ \pm\{6\}\right\}
\end{align*}
is a $B_{12}$-partition consisting of $3$ block pairs and the
zero-block $\{\pm 1,\,\pm 2,\,\pm 5,\,\pm 8,\,\pm 12\}$.

The total number of partitions of the set $[n]=\{1,2,\ldots,n\}$ is
called the Bell number and is denoted by $B_n$, see Rota
\cite{Rot64}. The type $B$ analogue of the Bell numbers are the
Dowling numbers $\left|\Pi_n^B\right|$, where $\Pi_n^B$ denotes the
set of $B_n$-partitions. The sequence $\{|\Pi_n^B|\}_{n\ge0}$ is
listed as A007405 in \cite{Slo}:
\[
1,\ 2,\ 6,\ 24,\ 116,\ 648,\ 4088,\ 28640,\ 219920,\ 1832224,\ \ldots.
\]

Let $\pi$ and $\pi'$ be two partitions of the set $[n]$. We say that
$\pi$ {\em refines} $\pi'$ if every block of $\pi$ is contained in
some block of $\pi'$. The refinement relation is a partial ordering
of the set $\Pi_n$ of all partitions of $[n]$. Define the {\em
meet}, denoted $\pi\wedge\pi'$, to be the largest partition which
refines both $\pi$ and $\pi'$. Define their {\em join}, denoted
$\pi\vee\pi'$, to be the smallest partition which is refined by both
$\pi$ and $\pi'$. The poset $\Pi_n$ is a lattice with the minimum
element $\hat0=\{\{1\},\{2\},\ldots,\{n\}\}$. We say that the
partitions $\pi_1,\pi_2,\ldots,\pi_r$ intersect minimally if
$\pi_1\wedge\pi_2\wedge\cdots\wedge\pi_r=\hat0$.

Pittel \cite{Pit00} has found a formula for the number of partitions minimally intersecting a given
 partition. He also computed the number of
minimally intersecting $r$-tuples of partitions.

\begin{thm}
Let $\pi$ be a partition of $[n]$, and let $i_1,\ldots,i_k$ be the sizes
of the blocks of $\pi$ listed in any order.
Then the number of partitions intersecting $\pi$ minimally equals
\[
N(\pi)=\bfi!\left[\bfx^\bfi\right]%
\exp\left(\prod_{\alpha\in[k]}(1+x_\alpha)-1\right),
\]
where $\bfi!=\prod_{\alpha\in[k]}i_\alpha!$ and $\left[\bfx^\bfi\right]$
stands for the coefficient of $\bfx^\bfi=\prod_{\alpha\in[k]}x_\alpha^{i_\alpha}$ of a
power series in $x_1, x_2, \ldots, x_k$.
Let $r\ge2$. The number $N_{n,r}$ of minimally intersecting $r$-tuples $(\pi_1,\pi_2,\ldots,\pi_r)$ of partitions is given by
\[
N_{n,r}=\frac{1}{e^r}\sum_{k_1,\ldots,k_r\ge0}%
\frac{(k_1k_2\cdots k_r)_n}{k_1!k_2!\cdots k_r!},
\]
where the notation $(x)_n=x(x-1)\cdots(x-n+1)$ denotes the falling factorial.
\end{thm}
By taking $\pi=\hat0$, the above formula reduces to Dobi\'{n}ski's
formula
\begin{equation}\label{Dobinski}
B_n=\frac{1}{e}\sum_{k\ge0}\frac{k^n}{k!},
\end{equation}
see Rota \cite{Rot64}. Wilf has obtained the following alternative
formula
\begin{equation}\label{Wilf}
N_{n,r}=\sum_{j=1}^nB_j^rs(n,j),
\end{equation}
where $s(n,j)$ is the Stirling number of the first kind.
Denote the generating function of $N_{n,r}$ by
\[
M_r(x)=\sum_{n\ge0}N_{n,r}\frac{x^n}{n!}.
\]
Canfield~\cite{Can01} has established the following connection
between $M_r(x)$ and the Bell numbers:
\begin{equation}\label{Canfield}
M_r\left(e^x-1\right)=\sum_{n\ge0}B_n^r\frac{x^n}{n!}.
\end{equation}
We shall give  type $B$ analogues of (\ref{Wilf}) and
(\ref{Canfield}) based on type $B$ partitions without zero block.

This paper is organized as follows. In Section 2, we give an expression
for the number of $B_n$-partitions that minimally intersect a
$B_n$-partition $\pi$ of a given type, which contains Benoumhani's
formula for the Dowling number as a special case. Moreover, we obtain a formula for the number of minimally
intersecting $r$-tuples of $B_n$-partitions. In Section 3, we consider
the enumeration of  minimally intersecting $r$-tuples of
$B_n$-partitions without zero-block, and give two formulas in analogy to
(\ref{Wilf}), and (\ref{Canfield}).


\section{Minimally intersecting $B_n$-partitions}

The main objective of this section is to derive a formula for the
number of minimally intersecting $r$-tuples of $B_n$-partitions.
If $\pi\in\Pi_n^B$ has a zero-block $Z=\{\pm i_1,\pm i_2,\ldots,\pm
i_k\}$, we say that $Z$ is of {\em half-size} $k$. The partition $
\hat0^B=\{\{1\},\,\{-1\},\,\{2\},\,\{-2\},\,\ldots,\,\{n\},\,\{-n\}\}
$ is called the minimal partition, and
$\hat1^B=\{\{\pm1,\,\pm2,\,\ldots,\,\pm n\}\}$ is called  the
maximal partition. We say that $\pi_1,\pi_2,\ldots,\pi_r$ are {\em
minimally intersecting} if
$\pi_1\wedge\pi_2\wedge\cdots\wedge\pi_r=\hat0^B$.

Let $\bfj=(j_1,j_2,\ldots,j_k)$ be a composition  of $n$. Let $\pi$
be a $B_n$-partition consisting of $k$ block pairs and a zero-block
of half-size $i_0$. For the purpose of enumeration, we often assume
that the block pairs of $\pi$ are ordered subject to certain
convention. We say that $\pi$ is of {\em type} $(i_0;\bfj)$ if the
block pairs of $\pi$ are ordered such that the $i$-th block pair is
of length $j_i$.

We first consider the problem of counting the number of
$B_n$-partitions with $l$ block pairs which minimally intersects a
given $B_n$-partition. As a special case, we are led to Benoumhani's
formula for the Dowling number
\begin{equation}\label{PinB}
\left|\Pi_n^B\right|=\frac{1}{\sqrt{e}}\sum_{k\ge0}\frac{(2k+1)^n}{(2k)!!},
\end{equation}
in analogy to Dobi\'{n}ski's formula~(\ref{Dobinski}). Next, we find
a formula for the number of ordered pairs of minimally intersecting
$B_n$-partitions. In general, we give a formula for the number of
minimally intersecting $r$-tuples of $B_n$-partitions.

\begin{thm}\label{thm_ParB}
Let $\pi$ be a $B_n$-partition consisting of a zero-block of
half-size $i_0$ (allowing $i_0=0$) and $k$ block pairs of sizes
$i_1,i_2,\ldots,i_k$ $(k\ge1)$ listed in any order. For any $l\ge1$,
the number of $B_n$-partitions $\pi'$ containing exactly $l$ block
pairs that intersect $\pi$ minimally equals
\begin{equation}\label{NBpil}
N^B(\pi;l)=\frac{\mathbf{i}!}{(2l-2i_0)!!}%
\sum_{\bfi'}\left[\mathbf{x^{i'}}\right]%
\left(\prod_{\alpha\in[k]}(1+x_\alpha)^2-1\right)^{l-i_0}%
\prod_{\alpha\in[k]}(1+x_\alpha)^{2i_0},
\end{equation}
where  $\bfi'$ runs over all vectors
$(i_1',i_2',\ldots,i_k')$ such that
$i_\alpha'\in\{i_\alpha,i_\alpha-1\}$ for any $\alpha\in[k]$, and
$\mathbf{x^{i'}}=\prod_{\alpha=1}^kx_\alpha^{i_\alpha'}$.
\end{thm}

For example,  $\Pi_2^B$ contains $6$ partitions:
\[
\hat0^B,\ \hat1^B,\ %
\{\pm\{1\},\{\pm2\}\},\ \{\pm\{2\},\{\pm1\}\},\ %
\{\pm\{1,2\}\},\ \{\pm\{1,-2\}\}.
\]
Let $\pi=\{\pm\{1\},\{\pm2\}\}$. We have $i_0=1$, $k=1$, and $i_1=1$.
For $l=1$, by~(\ref{NBpil}),
$N^B(\pi;1)=\sum_{i=0}^1\left[x^i\right](1+x)^{2}=3$. The three
$B_2$-partitions which contain exactly $1$ block pair and intersect
$\pi$ minimally are $\{\pm\{2\},\{\pm1\}\}$, $\{\pm\{1,2\}\}$, and
$\{\pm\{1,-2\}\}$.

\noindent{\it Proof of Theorem~\ref{thm_ParB}. }%
Let $Z_1$ be the zero-block of $\pi$, and  $\pm B_1,\pm
B_2,\ldots,\pm B_k$ be the block pairs of $\pi$. Let $Z_2$ be the
zero-block of $\pi'$, and $\pm B_1',\pm B_2',\ldots,\pm B_l'$ be the
block pairs of $\pi'$. To ensure that $\pi$ and $\pi'$ are minimally
intersected, it is necessary to characterize the intersecting
relations for all pairs $(B,B')$ where $B$ is a block of $\pi$ and
$B'$ is a block of $\pi'$.

First, we observe that the intersection $B\cap B'$ contains at most
one element subject to the minimally intersecting property. In
particular, $Z_1\cap Z_2=\emptyset$. If $B=Z_1$ and $B'\ne Z_2$,
then the two intersections $Z_1\cap B'$ and $Z_1\cap (-B')$ are a
pair of opposite subsets. This observation allows us to disregard
$Z_1\cap (-B')$ in our consideration. Since the cardinality of
$B\cap B'$ is either zero or one, we can represent $B\cap B'$ by
\[
F(k;l)\prod_{\beta\in[l]}(1+z_1w_\beta)%
\prod_{\alpha\in[k]}(1+x_\alpha z_2),
\]
where
\begin{equation}\label{Fkl}
F(k;l)=\prod_{\alpha\in[k],\,\beta\in[l]}%
(1+x_\alpha y_\beta)(1+x_\alpha \bar y_\beta).
\end{equation}
Here we use $x_i$ ($w_i$, resp.) to represent one of the two blocks
in the $i$-th block pair of $\pi$ ($\pi'$, resp.), and we use $y_i$ and $\bar
y_i$ to represent the two blocks in the $i$-th block pair of $\pi'$.

The above argument allows us to generate all $B_n$-partitions that
minimally meet with $\pi$. Let us consider the generating function
of such $B_n$-partitions. Set
\[
\begin{array}{lll}
\bfx=(x_1,x_2,\ldots,x_k),\quad%
&\bfi=(i_1,i_2,\ldots,i_k),\quad%
&\bfx^{\bfi}=\prod\limits_{\alpha\in[k]}x_\alpha^{i_\alpha};\\[16pt]
\bfw=(w_1,w_2,\ldots,w_l),\quad%
&\bfa=(a_1,a_2,\ldots,a_l),\quad%
&\bfw^{\bfa}=\prod\limits_{\beta\in[l]}w_\beta^{a_\beta};\\[16pt]
\bfy=(y_1,y_2,\ldots,y_l),\quad%
&\bfb=(b_1,b_2,\ldots,b_l),\quad%
&\bfy^{\bfb}=\prod\limits_{\beta\in[l]}y_\beta^{b_\beta};\\[16pt]
\bar\bfy=(\bar y_1,\bar y_2,\ldots,\bar y_l),\quad%
&\bfc=(c_1,c_2,\ldots,c_l),\quad%
&\bar\bfy^{\bfc}=\prod\limits_{\beta\in[l]}\bar y_\beta^{c_\beta}.
\end{array}
\]
Let $j_0$ be a nonnegative integer and $\bfj=(j_1,j_2,\ldots,j_l)$ a composition of
$n-j_0$. Denote by $N^B(\pi;j_0,\bfj)$ the number of $B_n$-partitions
$\pi'$ of type $(j_0;\bfj)$ such that $\pi'$ meets $\pi$ minimally.
In the above notation, we have
\begin{equation}\label{NBpij0j}
N^B(\pi;j_0,\bfj)=c\cdot\sum_{\bfa+\bfb+\bfc=\bfj}%
\left[\bfx^\bfi z_1^{i_0}z_2^{j_0}\bfw^\bfa\bfy^\bfb\bar{\bfy}^\bfc\right]%
F(k;l)\prod_{\beta\in[l]}(1+z_1w_\beta)%
\prod_{\alpha\in[k]}(1+x_\alpha z_2),
\end{equation}
where $c=\bfi!\cdot(2i_0)!!/(2l)!!$. Denote by ${S\choose m}$ the
collection of all $m$-subsets of $S$. Since
\begin{align}
\left[z_1^{i_0}\right]\prod_{\beta\in[l]}(1+z_1w_\beta)%
&=\sum_{Y\in{[l]\choose i_0}}\prod_{\beta\in Y}w_\beta,\label{z1}\\
\left[z_2^{j_0}\right]\prod_{\alpha\in[k]}(1+x_\alpha z_2)%
&=\sum_{X\in{[k]\choose j_0}}\prod_{\alpha\in X}x_\alpha,\label{z2}
\end{align}
substituting~(\ref{z1}) and~(\ref{z2}) into~(\ref{NBpij0j}), we
obtain that
\begin{align*}
N^B(\pi;j_0,\bfj)&=c\cdot\sum_{\bfa+\bfb+\bfc=\bfj}%
\left[\bfx^\bfi\bfw^\bfa\bfy^\bfb\bar{\bfy}^\bfc\right]%
\left(\sum_{Y\in{[l]\choose i_0}}\prod_{\beta\in Y}w_\beta\right)%
\left(\sum_{X\in{[k]\choose j_0}}\prod_{\alpha\in X}x_\alpha\right)F(k;l)\\
&=c\cdot\sum_{X,\,Y,\,\bfb}\left[\bfy^\bfb%
\prod_{\alpha\in[k]}x_\alpha^{i_\alpha-\chi(\alpha\in X)}%
\prod_{\beta\in[l]}\bar{y}_\beta^{j_\beta-b_\beta-\chi(\beta\in Y)}\right]F(k;l),
\end{align*}
where $\chi$ is the characteristic function defined by $\chi(P)=1$
if $P$ is true, and $\chi(P)=0$ otherwise. Therefore
\begin{equation}\label{eq9}
N^B(\pi;l)=\sum_{j_0+j_1+\cdots+j_l=n\atop{j_0\ge0,\,j_1,\ldots,j_l\ge1}}
N^B(\pi;j_0,\bfj)%
=c\cdot\sum_{j_0,\,X}%
\left[\prod_{\alpha}x_\alpha^{i_\alpha-\chi(\alpha\in X)}\right]%
\sum_{j_0+j_1+\cdots+j_l=n\atop{j_1,\ldots,j_l\ge1}}f(\bfj),
\end{equation}
where
\[
f(\bfj)=\sum_{Y,\,\bfb}\left[\bfy^\bfb\prod_{\beta}%
\bar y_\beta^{j_\beta-b_\beta-\chi(\beta\in Y)}\right]F(k;l).
\]
 In view of the expression~(\ref{Fkl}), the total degree of
$x_\alpha$'s agrees with the sum of the total degrees of
$y_\beta$'s and $\bar y_\beta$'s in $F(k;l)$. In other words,
\[
\sum_{\alpha\in[k]}i_\alpha-\chi(\alpha\in X)=\sum_{\beta\in[l]}b_\beta+(j_\beta-b_\beta-\chi(\beta\in
Y)),
\]
namely, $j_0+j_1+\cdots+j_l=n$. So we may drop this condition in the
inner summation of~(\ref{eq9}). For any $A\subseteq[l]$, let
\[
S(A)=\sum_{j_1,\ldots,j_l\ge0\atop{j_\beta=0\textrm{ if }\beta\not\in A}}f(\bfj)%
=\sum_{Y}\sum_{b_\gamma,j_\gamma\ge0\atop{\gamma\in A}}%
\left[\prod_{\gamma\in A}y_\gamma^{b_\gamma}%
\bar y_\gamma^{j_\gamma-b_\gamma-\chi(\gamma\in Y)}\right]F(k;A),
\]
where
\[
F(k;A)=\prod_{\alpha\in[k],\,\gamma\in A}%
(1+x_\alpha y_\gamma)(1+x_\alpha \bar y_\gamma).
\]
Since $j_\gamma$ and $b_\gamma$ run over all nonnegative
integers, the exponent $j_\gamma-b_\gamma-\chi(\gamma\in
Y)$ can considered as a summation index. It follows that
\begin{align*}
S(A)=\sum_{Y\in{A\choose i_0}}\sum_{b_\gamma,c_\gamma\ge0,\,\gamma\in A}%
\left[\prod_{\gamma\in A}y_\gamma^{b_\gamma}\bar y_\gamma^{c_\gamma}\right]%
F(k;A)={|A|\choose i_0}\prod_{\alpha\in[k]}(1+x_\alpha)^{2|A|}.
\end{align*}
By the principle of inclusion-exclusion, we have
\begin{align*}
\sum_{j_1,\ldots,j_l\ge1}f(\bfj)%
&=\sum_{A\subseteq[l]}(-1)^{l-|A|}S(A)%
=\sum_{m}{l\choose m}(-1)^{l-m}{m\choose i_0}%
\prod_{\alpha\in[k]}(1+x_\alpha)^{2m}\\
&={l\choose i_0}\prod_{\alpha\in[k]}(1+x_\alpha)^{2i_0}%
\left(\prod_{\alpha\in[k]}(1+x_\alpha)^2-1\right)^{l-i_0}.
\end{align*}
Now, employing (\ref{eq9}) we find that $N^B(\pi;l)$ equals
\begin{equation}\label{eq2}
\frac{\bfi!}{(2l-2i_0)!!}\sum_{X\subseteq[k]}%
\left[\prod_{\alpha\in[k]}x_\alpha^{i_\alpha-\chi(\alpha\in X)}\right]%
\prod_{\alpha\in[k]}(1+x_\alpha)^{2i_0}%
\left(\prod_{\alpha\in[k]}(1+x_\alpha)^2-1\right)^{l-i_0},
\end{equation}
which can be rewritten in the form of (\ref{NBpil}). This completes
the proof. \qed

The formula~(\ref{eq2}) will also be used in the proof of
Corollary~\ref{cor_NpilD}. Summing~(\ref{NBpil}) over $l\ge i_0$, we
obtain the following formula.

\begin{cor}
The number $N^B(\pi)$ of $B_n$-partitions that minimally intersect
$\pi$ is
\begin{equation}\label{NBpi}
N^B(\pi)=\frac{\bfi!}{\sqrt{e}}%
\sum_{\bfi'}\left[\mathbf{x^{i'}}\right]F(\bfx),
\end{equation}
where
\begin{equation}\label{Fx}
F(\bfx)=\left(\prod_{\alpha\in[k]}(1+x_\alpha)^{2i_0}\right)%
\exp\left(\frac{1}{2}\prod_{\alpha\in[k]}(1+x_\alpha)^2\right).
\end{equation}
\end{cor}

Setting $\pi=\hat0^B$, (\ref{NBpi}) reduces to (\ref{PinB}), since
\[ N^B(\hat0^B)%
=\frac{1}{\sqrt{e}}\sum_{i_\alpha'\in\{0,1\}}%
\left[x_1^{i_1'}\cdots x_n^{i_n'}\right]%
\sum_{j\ge0}\frac{1}{(2j)!!}\prod_{\alpha=1}^n(1+x_\alpha)^{2j}.
\]

In fact, the number $N^B(\pi)$ can be expressed in terms of an infinite sum.

\begin{cor}
\begin{equation}\label{nb2}
N^B(\pi)=\frac{1}{\sqrt{e}}\sum_{j\ge0}\frac{(2i_0+2j+1)!^k}{(2j)!!}%
\prod_{\alpha\in[k]}\frac{1}{(2i_0+2j+1-i_\alpha)!}.
\end{equation}
\end{cor}

\pf  From (\ref{Fx}) it follows that
\[
F(x)=\sum_{j\ge0}\frac{1}{(2j)!!}\prod_{\alpha\in[k]}(1+x_\alpha)^{2(i_0+j)}.
\]
Hence
\begin{align*}
N^B(\pi)&=\frac{\bfi!}{\sqrt{e}}\sum_{j\ge0}\frac{1}{(2j)!!}%
\prod_{\alpha\in[k]}\left({2(i_0+j)\choose i_\alpha}+{2(i_0+j)\choose i_\alpha-1}\right)\\
&=\frac{\bfi!}{\sqrt{e}}\sum_{j\ge0}\frac{1}{(2j)!!}%
\prod_{\alpha\in[k]}{2(i_0+j)+1\choose i_\alpha},
\end{align*}
which gives (\ref{nb2}).  This completes the proof. \qed

\begin{cor}\label{cor_Nn2Bi0k}
Let $N_{n,2}^B(i_0;k)$ denote the number of ordered pairs $(\pi,\pi')$ of minimally
intersecting $B_n$-partitions such that $\pi$ consists of exactly $k$ block pairs and a zero-block of half-size $i_0$ (allowing $i_0=0$). Then
\begin{equation}\label{Nn2Bi0k}
N_{n,2}^B(i_0;k)=\frac{(2n)!!}{(2i_0)!!(2k)!!\sqrt{e}}%
\left[x^{n-i_0}\right]\sum_{j\ge0}\frac{1}{(2j)!!}%
\left((1+x)^{2i_0+2j+1}-1\right)^k.
\end{equation}
\end{cor}

\pf By a simple combinatorial argument we see that the number of
$B_n$-partitions of type $(i_0;i_1,\ldots,i_k)$ equals
\[
c={n\choose i_0,i_1,\ldots,i_k}\frac{2^{n-i_0-k}}{k!}%
=\frac{(2n)!!}{(2i_0)!!(2k)!!}\cdot\frac{1}{\mathbf{i}!}.
\]
Thus by~(\ref{NBpi}), we have
\begin{equation}\label{eq1}
N_{n,2}^B(k)%
=\sum_{i_0+i_1+\cdots+i_k=n\atop{i_1,\ldots,i_k\ge1}}c\cdot N^B(\pi)%
=\frac{(2n)!!}{(2i_0)!!(2k)!!\sqrt{e}}%
\sum_{i_0+i_1+\cdots+i_k=n\atop{i_1,\ldots,i_k\ge1}}%
\sum_{\bfi'}\left[\mathbf{x^{i'}}\right]F(\bfx).
\end{equation}
For any $A\subseteq[k]$, consider
\begin{align*}
S(A)=\sum_{i_0+i_1+\cdots+i_k=n%
\atop{i_1,\ldots,i_k\ge0\atop{i_\alpha=0\textrm{ if }\alpha\not\in A}}}%
\sum_{\bfi'}\left[\mathbf{x^{i'}}\right]F(\bfx)%
=\sum_{i_0+\sum_{\alpha\in A}i_\alpha=n%
\atop{i_\alpha\ge0,\,\alpha\in A}}%
\sum_{\bfi'|_A}\left[\bfx^{\bfi'}\big|_A\right]F\left(\bfx\big|_A\right),
\end{align*}
where $\bfx\big|_A$ (resp. $\bfi'|_A$) denotes the vector obtained
by removing all $x_\alpha$ (resp. $i_\alpha'$) such that
$\alpha\not\in A$ from the vector $\bfx$ (resp. $\bfi'$). Let $t$ be
the number of $\alpha$'s such that $i_\alpha'=i_\alpha-1$ in the
inner summation. Noting that
\[
F\left(\bfx\big|_A\right)=\left(\prod_{\alpha\in A}(1+x_\alpha)^{2i_0}\right)%
\exp\left(\frac{1}{2}\prod_{\alpha\in A}(1+x_\alpha)^2\right),
\]
we can transform $S(A)$ to
\begin{align*}
S(A)&=\left(\sum_{t}{|A|\choose t}\left[x^{n-i_0-t}\right]\right)%
(1+x)^{2i_0|A|}\exp\left(\frac{1}{2}(1+x)^{2|A|}\right)\\
&=\left[x^{n-i_0}\right]%
(1+x)^{(2i_0+1)|A|}\exp\left(\frac{1}{2}(1+x)^{2|A|}\right).
\end{align*}
In view of the principle of inclusion-exclusion, we deduce
from~(\ref{eq1}) that
\[
N_{n,2}^B(k)=\frac{(2n)!!}{(2i_0)!!(2k)!!\sqrt{e}}%
\sum_{A\subseteq[k]}(-1)^{k-|A|}S(A),
\]
which gives~(\ref{Nn2Bi0k}). This completes the proof. \qed

Summing over $0\le k\le n-i_0$ and $0\le i_0\le n$, we obtain the
number of ordered pairs of minimally intersecting $B_n$-partitions.

\begin{cor}
The number $N_{n,2}^B$ of ordered pairs $(\pi,\pi')$ of minimally
intersecting $B_n$-partitions is given by
\[
N_{n,2}^B=\frac{2^n}{e}\sum_{k,l\ge0}\frac{(2kl+k+l)_n}{(2k)!!(2l)!!}.
\]
\end{cor}

For example, $N_{1,2}^B=3$, $N_{2,2}^B=23$, $N_{3,2}^B=329$,
$N_{4,2}=6737$. In general, we have the following theorem, which is the main result of this paper.

\begin{thm}\label{thm_main}
Let $r\ge2$. The number of minimally intersecting $r$-tuples $(\pi_1,\pi_2,\ldots,\pi_r)$
of $B_n$-partitions equals
\begin{equation}\label{NnrB}
N_{n,r}^B=\frac{2^n}{e^{r/2}}\sum_{l_1,l_2,\ldots,l_r}%
\frac{(f_r)_n}{(2l_1)!!(2l_2)!!\cdots(2l_r)!!},
\end{equation}
where
\[
f_r=\frac{1}{2}\left(\prod_{t\in[r]}(2l_t+1)-1\right).
\]
\end{thm}

\pf For any $t\in[r]$, let $i_t$ be an nonnegative integer and
$\bfj_t=(j_{t,1},\,j_{t,2},\,\ldots,\,j_{t,k_t})$ be a composition
of $n$. Let $\pi_t$ be a $B_n$-partition of type $(i_t;\bfj_t)$. The
condition that $\pi_1,\pi_2,\ldots,\pi_r$ are minimally intersecting
leads us to consider the intersecting relations for all $r$-tuples
$(B_1,B_2,\ldots,B_r)$ where $B_t$ is a block of $\pi_t$.

First, we observe that the intersection
\begin{equation}\label{intersection}
B_1\cap B_2\cap\cdots\cap B_r
\end{equation}
contains at most one element because of the minimally intersecting
requirement. In particular, (\ref{intersection}) is empty when
$B_1,B_2,\ldots,B_r$ are all zero-blocks. We now consider the case
that not all of $B_1, B_2, \ldots, B_r$ are zero-blocks. In this
case, there exists a number $t\in[r]$ such that $B_1,\ldots,B_{t-1}$
are zero-blocks but $B_t$ is a non-zero-block. This number $t$ will
play a key role in determining the
intersection~(\ref{intersection}).

In fact, the partial intersection $B_1\cap B_2\cap\cdots\cap
B_{t-1}$ is of the form $\{\pm i_1,\ldots,\pm i_j\}$. Thus for any
non-zero-block $B$ of $\pi_t$, the two intersections
\[
B_1\cap\cdots\cap B_{t-1}\cap B\quad \mbox{and}\quad
B_1\cap\cdots\cap B_{t-1}\cap(-B)
\]
form a pair of opposite subsets. This observation allows us
to consider $B$ as a representative of the block pair $\pm B$.
Since the cardinality of the intersection~(\ref{intersection}) is
either zero or one, we can represent~(\ref{intersection}) by
\begin{equation}\label{f}
f=1+z_1\cdots z_{t-1}x_{t,\alpha_t}Y_{t+1}\cdots Y_r,
\end{equation}
where
\[
Y_p\in\left\{z_p,\ y_{p,1},\ \bar y_{p,1},\ \ldots,\ y_{p,k_p},\
\bar y_{p,k_p}\right\}
\]
for $p\ge t+1$. Here we use $z_i$ to represent the zero-block of
$\pi_i$, $x_{t,i}$ to represent one of the two blocks in the
$i$-th block pair of $\pi_t$, $y_{p,i}$ and $\bar y_{p,i}$ to
represent the two blocks in the $i$-th block pair of $\pi_p$. Let
\[
\begin{array}{lll}
\bfx_t=(x_{t,1},\,\ldots,\,x_{t,k_t}),\quad%
&\bfa_t=(a_{t,1},\,\ldots,\,a_{t,k_t}),\quad%
&\bfx_s^{\bfa_s}=\prod\limits_{i\in[k_s]}x_{s,i}^{a_{s,i}};\\[16pt]
\bfy_t=(y_{t,1},\,\ldots,\,y_{t,k_t}),\quad%
&\bfb_t=(b_{t,1},\,\ldots,\,b_{t,k_t}),\quad%
&\bfy_s^{\bfb_s}=\prod\limits_{i\in[k_s]}y_{s,i}^{b_{s,i}};\\[16pt]
\bar{\bfy}_t=(\bar y_{t,1},\,\ldots,\,\bar y_{t,k_t}),\quad%
&\bfc_t=(c_{t,1},\,\ldots,\,c_{t,k_t}),\quad%
&\bar\bfy_s^{\bfc_s}=\prod\limits_{i\in[k_s]}\bar y_{s,i}^{c_{s,i}}.
\end{array}
\]
Denote by $N^B(\pi_1;i_2,\bfj_2;\ldots;i_r,\bfj_r)$  the number of
$(r-1)$-tuples $(\pi_2,\ldots,\pi_r)$ of $B_n$-partitions such that
$\pi_s$ ($2\le s\le r$) is of type $(i_s,\bfj_s)$ and
$\pi_1,\pi_2,\ldots,\pi_r$ intersect minimally. In the notation of
$f$ in~(\ref{f}), we get
\[
N^B(\pi_1;i_2,\bfj_2;\ldots;i_r,\bfj_r)%
=c\left[\bfx_1^{\bfj_1}z_1^{i_1}\right]%
\sum_{\bfa_s+\bfb_s+\bfc_s=\bfj_s\atop{2\le s\le r}}%
\left[\bfx_s^{\bfa_s}\bfy_s^{\bfb_s}\bar\bfy_s^{\bfc_s}z_s^{i_s}\right]F_r%
\]
where
\begin{align*}
c&=\bfj_1!\cdot(2i_1)!!\prod_{2\le s\le r}(2k_s)!!^{-1},\\
F_r&=\prod_{t\in[r]}\prod_{\alpha_t\in\left[k_t\right]}%
\prod_{Y_p\in\left\{z_p,y_{p,1},\bar y_{p,1},\ldots,y_{p,k_p},\bar y_{p,k_p}\right\}%
\atop{t+1\le p\le r}}f.
\end{align*}

Now, let $N^B(\pi_1,k_2,\ldots,k_r)$ be the number of $(r-1)$-tuples $(\pi_2,\ldots,\pi_r)$ of $B_n$-partitions such that $\pi_s$ contains exactly $k_s$ block pairs and $\pi_1,\pi_2,\ldots,\pi_r$ intersect minimally. Then
\begin{equation}\label{eq10}
N^B(\pi_1,k_2,\ldots,k_r)%
=\sum_{i_s\ge0,\,j_{s,1},\ldots,j_{s,k_s}\ge1%
\atop{j_{s,1}+\cdots+j_{s,k_s}+i_s=n}}%
N^B(\pi_1;i_2,\bfj_2;\ldots;i_r,\bfj_r)
\end{equation}
We claim that the condition $j_{s,1}+\cdots+j_{s,k_s}+i_s=n$ can be
dropped in the above summation. In fact, the factor $f$ in~(\ref{f})
contributes to $x_1$ or $z_1$ at most once with respect to the
degree, and the contribution of $f$ to $x_1$ or $z_1$ equals the
contribution of $f$ to $\bfx_s$, $\bfy_s$, $\bar\bfy_s$, or $z_s$,
for any $2\le s\le r$. Therefore the sum of the degrees of $\bfx_s$,
$\bfy_s$, $\bar\bfy_s$, and $z_s$, equals the sum of the degrees of
$x_1$ and $z_1$, that is, for any $2\le s\le r$,
\begin{equation}\label{condition}
i_s+j_{s,1}+\cdots+j_{s,k_s}=i_1+j_{1,1}+\cdots+j_{1,k_1}=n
\end{equation}
Hence we can ignore the conditions~(\ref{condition})
in~(\ref{eq10}). This implies that
\[
N^B(\pi_1,k_2,\ldots,k_r)%
=c\left[\bfx_1^{\bfj_1}z_1^{i_1}\right]%
\sum_{i_s\ge0,\,\bfa_s+\bfb_s+\bfc_s\ge{\bf1}}%
\left[\bfx_s^{\bfa_s}\bfy_s^{\bfb_s}\bar\bfy_s^{\bfc_s}z_s^{i_s}\right]F_r.
\]
where $\bfa_s+\bfb_s+\bfc_s\ge{\bf1}$ indicates that
$a_{s,h_s}+b_{s,h_s}+c_{s,h_s}\ge1$ for any $1\le h_s\le k_s$. We
will compute  $\sum%
\left[\bfx_s^{\bfa_s}\bfy_s^{\bfb_s}\bar\bfy_s^{\bfc_s}z_s^{i_s}\right]F_r$
for $s=2,3,\ldots,r$ by the following procedure. First,  for $s=2$, we have
\[
\sum_{i_2\ge0,\,\bfa_2+\bfb_2+\bfc_2\ge\bf1}%
\left[\bfx_2^{\bfa_2}\bfy_2^{\bfb_2}\bar{\bfy}_2^{\bfc_2}z_2^{i_2}\right]F_r
=\sum_{l_2}{k_2\choose l_2}(-1)^{k_2-l_2}F_{r,2},
\]
where
\[
F_{r,2}
=\prod_{\alpha_1,Y_p}(1+x_1^{\alpha_1}Y_3\cdots Y_r)^{2l_2+1}%
\prod_{Y_p}(1+z_1Y_3\cdots Y_r)^{l_2}
\prod_{t\ge3,\,\alpha_t,\,Y_p}%
(1+z_1z_3\cdots z_{t-1}x_t^{\alpha_t}Y_{t+1}\cdots Y_r).
\]
So $N^B(\pi_1,k_2,\ldots,k_r)$ equals
\begin{equation}\label{eq3}
c\left[\bfx_1^{\bfj_1}z_1^{i_1}\right]%
\sum_{l_2}{k_2\choose l_2}(-1)^{k_2-l_2}%
\sum_{i_s\ge0,\,\bfa_s+\bfb_s+\bfc_s\ge{\bf1}\atop{3\le s\le r}}%
\left[\bfx_s^{\bfa_s}\bfy_s^{\bfb_s}\bar\bfy_s^{\bfc_s}z_s^{i_s}\right]%
F_{r,2}.
\end{equation}
To  compute the inner summation,   let
\[
g_s=\frac{1}{2}\left(\prod_{2\le i\le s}(2l_i+1)-1\right).
\]
For any $s\ge2$, it is clear that
\[
(2l_{s+1}+1)g_s+l_{s+1}=g_{s+1}.
\]
Starting with (\ref{eq3}),  we can continue the above procedure to deduce that for $2\le h\le r-1$,
\[
N^B(\pi_1,k_2,\ldots,k_r)%
=c\left[\bfx_1^{\bfj_1}z_1^{i_1}\right]%
\sum_{l_2,\ldots,l_h}%
\prod_{2\le i\le h}{k_i\choose l_i}(-1)^{k_i-l_i}%
\sum_{i_s\ge0,\,\bfa_s+\bfb_s+\bfc_s\ge{\bf1}\atop{h+1\le s\le r}}%
\left[\bfx_s^{\bfa_s}\bfy_s^{\bfb_s}\bar\bfy_s^{\bfc_s}z_s^{i_s}\right]%
F_{r,h},
\]
where
\begin{align*}
F_{r,h}%
=&\prod_{\alpha_1,Y_p}%
\left(1+x_1^{\alpha_1}Y_{h+1}\cdots Y_r\right)^{\prod_{2\le i\le
h}(2l_i+1)}%
\prod_{Y_p}\left(1+z_1Y_{h+1}\cdots Y_r\right)^{g_h}\\
&\cdot\prod_{t\ge h+1,\,\alpha_t,\,Y_p}%
\left(1+z_1z_{h+1}\cdots z_{t-1}x_t^{\alpha_t}Y_{t+1}\cdots
Y_r\right).
\end{align*}
In particular, for $h=r-1$, we have
\begin{equation}\label{r-1}
N^B(\pi_1,k_2,\ldots,k_r)%
=c\left[\bfx_1^{\bfj_1}z_1^{i_1}\right]%
\sum_{l_2,\ldots,l_{r-1}}%
\left(\prod_{2\le i\le r-1}{k_i\choose l_i}(-1)^{k_i-l_i}\right)G%
\end{equation}
where
\begin{align*}
G&=\sum_{\bfa_r+\bfb_r+\bfc_r\ge{\bf1}}%
\left[\bfx_r^{\bfa_r}\bfy_r^{\bfb_r}\bar\bfy_r^{\bfc_r}\right]%
\prod_{\alpha_1,Y_p}%
\left(1+x_1^{\alpha_1}\right)^{\prod_{2\le i\le r-1}(2l_i+1)}%
\prod_{Y_p}\left(1+z_1\right)^{g_{r-1}}%
\prod_{\alpha_r}\left(1+z_1x_r^{\alpha_r}\right)\\
&=\sum_{l_r}{k_r\choose l_r}(-1)^{k_r-l_r}(1+z_1)^{g_r}%
\prod_{\alpha_1}\left(1+x_1^{\alpha_1}\right)^{\prod_{2\le i\le
r}(2l_i+1)}.
\end{align*}
Since the number of $B_n$-partitions of type $\bfj_1$ equals
\[
c'={n\choose i_1}{n-i_1\choose\bfj_1}\frac{2^{n-i_1-k_1}}{k_1!}%
=\frac{(2n)!!}{(2i_1)!!(2k_1)!!\bfj_1!},
\]
by (\ref{r-1}), we obtain
\begin{align}
N_{n,r}^B&=\sum_{j_{1,1},\ldots,j_{1,k_1}\ge1%
\atop{i_1+j_{1,1}+\cdots+j_{1,k_1}=n}}c'%
\sum_{k_2,\ldots,k_r}N^B(\pi_1,k_2,\ldots,k_r)\notag\\
&=(2n)!!\sum_{k_2,\ldots,k_r\atop{l_2,\ldots,l_r}}%
\left(\prod_{2\le s\le r}{k_s\choose l_s}%
\frac{(-1)^{k_s-l_s}}{(2k_s)!!}\right)%
\sum_{i_1,k_1}\frac{1}{(2k_1)!!}\left[z_1^{i_1}\right](1+z_1)^{g_r}H\label{eq8}
\end{align}
where
\begin{align*}
H&=\sum_{i_1+j_{1,1}+\cdots+j_{1,k_1}=n%
\atop{j_{1,1},\,j_{1,2},\,\ldots,\,j_{1,k_1}\ge1%
}}\left[\bfx_1^{\bfj_1}\right]%
\prod_{\alpha_1}\left(1+x_1^{\alpha_1}\right)^{\prod_{2\le i\le
r}(2l_i+1)}\\
&=\sum_{l_1}{k_1\choose l_1}(-1)^{k_1-l_1}%
\left[x^{n-i_1}\right](1+x)^{l_1\prod_{2\le i\le r}(2l_i+1)}.
\end{align*}
Using the identity
\begin{equation}\label{identity}
\sum_{k}{k\choose l}\frac{(-1)^{k-l}}{(2k)!!}%
=\frac{e^{-1/2}}{(2l)!!},
\end{equation}
 we can simplify the summation over $k_1, k_2,\ldots,k_r\ge0$
in~(\ref{eq8}) to deduce that
\begin{align}
N_{n,r}^B&=(2n)!!\sum_{k_1,k_2,\ldots,k_r\atop{l_1,l_2,\ldots,l_r}}%
\left(\prod_{t\in[r]}{k_t\choose l_t}%
\frac{(-1)^{k_t-l_t}}{(2k_t)!!}\right)%
\sum_{i_1}\left[x^{n-i_1}z_1^{i_1}\right](1+z_1)^{g_r}%
(1+x)^{l_1\prod_{2\le i\le r}(2l_i+1)}\notag\\
&=\frac{(2n)!!}{e^{r/2}}\sum_{l_1,l_2,\ldots,l_r}%
\frac{1}{(2l_1)!!(2l_2)!!\cdots(2l_r)!!}\left[x^n\right](1+x)^{g_r+l_1\prod_{2\le
i\le r}(2l_i+1)}.\label{eq7}
\end{align}
To further simplify the above summation, we observe that
\begin{equation}\label{f_r}
g_r+l_1\prod_{2\le i\le
r}(2l_i+1)=\frac{1}{2}\left(\prod_{t\in[r]}(2l_t+1)-1\right).
\end{equation}
Substituting~(\ref{f_r}) into~(\ref{eq7}), we arrive at
(\ref{NnrB}). This completes the proof. \qed

For example, we have $N_{1,r}=2^r-1$ and $N_{2,3}^B=187$.

\section{Minimally intersecting $B_n$-partitions without zero-block}

In this section, we investigate the meet-semilattice of
$B_n$-partitions without zero-block. Note that the minimal
$B_n$-partition without zero-block is $\hat0^B$. Inspecting the
proof of Theorem~\ref{thm_ParB}, we can restrict our attention to
the set of $B_n$-partitions without zero-block  by setting $i_0=0$
and $X=\emptyset$ in~(\ref{eq2}).

\begin{cor}\label{cor_NpilD}
Let $\pi$ be a $B_n$-partition consisting of $k$ block pairs of
sizes $i_1,i_2,\ldots,i_k$ listed in any order. For a given $l\ge1$,
the number $N^D(\pi;l)$ of $B_n$-partitions $\pi'$ consisting of $l$
block pairs, which intersects $\pi$ minimally, is equal to
\begin{equation}\label{NpilD}
N^D(\pi;l)=\frac{\mathbf{i}!}{(2l)!!}\left[\mathbf{x^i}\right]%
\left(\prod_{\alpha\in[k]}(1+x_\alpha)^2-1\right)^l.
\end{equation}
The total number of $B_n$-partitions without zero-block that intersect $\pi$ minimally is given by
\begin{equation}\label{NDpi}
N^D(\pi)=\frac{\mathbf{i}!}{\sqrt{e}}\left[\bfx^\bfi\right]%
\exp\left(\frac{1}{2}\prod_{\alpha\in[k]}(1+x_\alpha)^2\right).
\end{equation}
\end{cor}

For example, let $n=3$, $\pi=\{\pm\{2\},\,\pm\{1,-3\}\}$ and $l=2$.
Then (\ref{NpilD}) yields $N^D(\pi;2)=5$. In fact, the
$B_n$-partitions consisting of $2$ block pairs which intersect
$\pi$ minimally are exactly the $5$ partitions consisting of two
block pairs except for $\pi$ itself.

Let $N_n$ be the number of $B_n$-partitions without zero-block.
Taking $\pi=\hat0^B$ in~(\ref{NDpi}), we obtain the following
formula.

\begin{cor} We have
\begin{equation}\label{N_n}
N_n=\frac{1}{\sqrt{e}}\sum_{k\ge0}\frac{(2k)^n}{(2k)!!}.
\end{equation}
\end{cor}

Let $N_n(k)$ denote the number of $B_n$-partitions containing $k$
block pairs but no zero-block. It should be noted that the formula
(\ref{N_n}) can be easily deduced from the relation
\[
N_n(k)=2^{n-k}S(n,k),
\]
where $S(n,k)$ are the Stirling numbers of the second kind, and the
following identity on the Bell polynomials \cite{Rio80, Rom84}:
\[
\sum_{k=0}^nS(n,k)x^k=\frac{1}{e^x}\sum_{k\ge0}\frac{k^n}{k!}x^k.
\]

The sequence $\{N_n\}_{n\ge0}$ is A004211 in \cite{Slo}:
\[
1,\ 1,\ 3,\ 11,\ 49,\ 257,\ 1539,\ 10299,\ 75905,\ 609441,\ \ldots.
\]

The proof of Corollary~\ref{cor_Nn2Bi0k} implies the following
corollary.

\begin{cor}
Let $N_{n,2}^D(k)$ denote the number of ordered pairs $(\pi,\pi')$ of minimally
intersecting $B_n$-partitions without zero-block such that $\pi$ consists of exactly $k$ block pairs. Then
\[
N_{n,2}^D(k)=\frac{(2n)!!}{(2k)!!\sqrt{e}}\left[x^n\right]%
\sum_{j\ge0}\frac{1}{(2j)!!}\left[(1+x)^{2j}-1\right]^k.
\]
The number $N_{n,2}^D$ of ordered pairs $(\pi,\pi')$ of
minimally intersecting $B_n$-partitions without zero-block is given by
\[
N_{n,2}^D=\frac{2^n}{e}\sum_{k,\,l\ge0}\frac{(2kl)_n}{(2k)!!\,(2l)!!}.
\]
\end{cor}

For example, $N_{1,2}^D=1$, $N_{2,2}^D=7$, $N_{3,2}^D=75$. The
following theorem is an analogue of Theorem~\ref{thm_main} with
respect to the meet-semilattice of $B_n$-partitions without
zero-block.

\begin{thm}
For $r\ge2$, the number of minimally intersecting $r$-tuples $(\pi_1,\pi_2,\ldots,\pi_r)$
of $B_n$-partitions without zero-block equals
\begin{equation}\label{NnrD}
N_{n,r}^D=\frac{2^n}{e^{r/2}}\sum_{k_1,k_2,\ldots,k_r}%
\frac{\left(2^{r-1}k_1k_2\cdots k_r\right)_n}{(2k_1)!!(2k_2)!!\cdots(2k_r)!!}.
\end{equation}
\end{thm}

\pf Let $1\le t\le r$. Let
$\bfj_t=(j_{t,1},j_{t,2},\ldots,j_{t,k_t})$ be a composition of $n$.
Assume that $\pi_t$ is of type $(0;\bfj_t)$. Let
$N^D(\pi_1,\bfj_2,\ldots,\bfj_r)$ be the number of $(r-1)$-tuples
$(\pi_2,\ldots,\pi_r)$ of such $B_n$-partitions such that $(\pi_1,
\pi_2, \ldots, \pi_r)$ is minimally intersecting. By the argument in
the proof of Theorem~\ref{thm_ParB}, we find
\begin{equation}\label{nd}
N^D(\pi_1,\bfj_2,\ldots,\bfj_r)%
=c\cdot\left[\bfx^{\bfj_1}\right]\sum_{\bfb_s+\bfc_s=\bfj_s}%
\left[\bfy_2^{\bfb_2}\bar{\bfy}_2^{\bfc_2}\cdots%
\bfy_r^{\bfb_r}\bar{\bfy}_r^{\bfc_r}\right]f(\bfj),
\end{equation}
where
\begin{align*}
c&=\bfj_1!\prod_{2\le s\le r}(2k_s)!!^{-1},\\
f(\bfj)&=\prod_{\alpha\in[k_1]%
\atop{Y_s\in\left\{y_{s,1},\,\bar y_{s,1},\,\ldots,\,y_{s,k_s},\,\bar y_{s,k_s}\right\}}}%
\left(1+x_\alpha Y_2Y_3\cdots Y_r\right).
\end{align*}

Let $N^D(\pi_1,k_2,\ldots,k_r)$ be the number of $(r-1)$-tuples
$(\pi_2,\ldots,\pi_r)$ of $B_n$-partitions such that $\pi_s$
consists of $k_s$ block pairs, and $\pi_1, \pi_2, \ldots, \pi_r$ are
minimally intersecting. It follows from (\ref{nd}) that
\begin{align*}
N^D(\pi_1,k_2,\ldots,k_r)%
&=c\cdot\left[\bfx^{\bfj_1}\right]%
\sum_{\bfb_s+\bfc_s=\bfj_s\ge{\bf1}}%
\left[\bfy_2^{\bfb_2}\cdots\bar{\bfy}_r^{\bfc_r}\right]f(\bfj)\\
&=\bfj_1!\sum_{l_2,\ldots,l_r}%
\left(\left[\bfx^{\bfj_1}\right]%
\prod_{\alpha\in[k_1]}(1+x_\alpha)^{2^{r-1}l_2\cdots l_r}\right)%
\prod_{2\le s\le r}{k_s\choose l_s}\frac{(-1)^{k_s-l_s}}{(2k_s)!!}.
\end{align*}
Consequently,
\begin{align*}
N_{n,r}^D&=\sum_{k_1}\frac{1}{(2k_1)!!}%
\sum_{j_{1,1}+\cdots+j_{1,k_1}=n\atop{j_{1,1},\ldots,j_{1,k_1}\ge1}}\frac{2^{n}n!}{\bfj_1!}%
\sum_{k_2,\ldots,k_r}N^D(\pi_1,k_2,\ldots,k_r)\\
&=(2n)!!\sum_{k_1,k_2,\ldots,k_r\atop{l_1,l_2,\ldots,l_r}}%
\prod_{1\le s\le r}{k_s\choose l_s}\frac{(-1)^{k_s-l_s}}{(2k_s)!!}%
[x^n](1+x)^{2^{r-1}l_1l_2\cdots l_r}.
\end{align*}%
Applying (\ref{identity}), we can restate the above formula in the form of (\ref{NnrD}).
This completes the proof. \qed

For example, when $n=2$ and $r=3$, by~(\ref{NnrD}) we find that $
N_{2,3}^D=25$. In fact, there are $3$ $B_2$-partitions without
zero-block, that is,
\[
0^B,\ \pi_1=\{\pm\{1,2\}\},\ \pi_2=\{\pm\{1,-2\}\}.
\]
Among all $3^3=27$ $3$-tuples of $B_2$-partitions without
zero-block, only $(\pi_1,\pi_1,\pi_1)$ and $(\pi_2,\pi_2,\pi_2)$ are
not minimally intersecting.

\begin{cor}
We have
\begin{equation}\label{WilfD}
N_{n,r}^D=\sum_{j=1}^nN_j^r 2^{n-j}s(n,j),
\end{equation}
where $s(n,j)$ are the Stirling numbers of the first kind. Moreover,
\begin{equation}\label{CanfieldD}
M_r^D\left(\frac{e^{2x}-1}{2}\right)=\sum_{n\ge0}N_n^r\frac{x^n}{n!},
\end{equation}
where
\[
M_r^D(x)=\sum_{n\ge0}N_{n,r}^D\frac{x^n}{n!}.
\]
\end{cor}

The formula~(\ref{WilfD}) can be considered as a type $B$ analogue
of Wilf's formula~(\ref{Wilf}), whereas (\ref{CanfieldD}) is
analogous to Canfield's formula~(\ref{Canfield}).

\noindent{\bf Acknowledgments.} This work was supported by the 973
Project, the PCSIRT Project of the Ministry of Education, and the
National Science Foundation of China.


\begin{thebibliography}{99}
\small \setlength{\itemsep}{-.8mm}

\bibitem{Ben96}%
M. Benoumhani, On Whitney numbers of Dowling lattices, Discrete
Math.,  159 (1996) 13--33.

\bibitem{BB05}%
A. Bj\"{o}rner and F. Brenti, Combinatorics of Coxeter Groups, 2005,
Springer Science+Business Media, Inc.

\bibitem{BW04}%
A. Bj\"{o}rner and M.L. Wachs, Geometrically constructed bases for
homology of partitions lattices of types $A$, $B$ and $D$, Electron.
J. Combin., 11 (2004) \#R3.

\bibitem{Can01}%
E.R. Canfield, Meet and join within the lattice of set partitions,
Electron. J. Combin., 8 (2001) \#R15.

\bibitem{Dow73}%
T.A. Dowling, A class of geometric lattices based on finite groups,
J. Combin. Theory Ser. B, 14 (1973) 61--86.

\bibitem{Hum90}%
J.E. Humphreys, Reflection Groups and Coxeter Groups, Cambridge
Studies in Advanced Mathematics 29, Cambridge Univ. Press,
Cambridge, 1990.

\bibitem{Pit00}%
B. Pittel, Where the typical set partitions meet and join, Electron.
J. Combin., 7 (2000) \#R5.

\bibitem{Rei97}%
V. Reiner, Non-crossing partitions for classical reflection groups,
Discrete Math., 177 (1997) 195--222.

\bibitem{Rio80}%
J. Riordan, An Introduction to Combinatorial Analysis, Wiley, New
York, 1980.

\bibitem{Rom84}%
S. Roman, The Umbral Calculus,  Academic Press, New York, 1984.

\bibitem{Rot64}%
G.C. Rota, The number of partitions of a set, Amer. Math. Monthly,
71 (1964) 498--504.

\bibitem{Slo}%
N.J.A. Sloane, The On-Line Encyclopedia of Integer Sequences,\\
{\tt http://www.research.att.com/\char 126 njas/sequences/}.

\end{thebibliography}
\end{document}